\documentclass[11pt,a4paper]{article}
\usepackage[T1]{fontenc}
\usepackage{lmodern,microtype}
\usepackage[margin=27mm]{geometry}
\usepackage{amsmath,amssymb,amsthm,mathtools,booktabs,longtable,array}
\usepackage{xcolor}
\usepackage[colorlinks=true,linkcolor=blue!45!black,citecolor=blue!45!black,urlcolor=blue!45!black]{hyperref}
\usepackage{bookmark}
\hypersetup{pdftitle={Euclidean View-Obstruction in Dimensions Eight to Ten and the Asymptotic Deficit},pdfauthor={Hu Tan and Ying Zhang}}
\newtheorem{theorem}{Theorem}[section]
\newtheorem{proposition}[theorem]{Proposition}
\newtheorem{lemma}[theorem]{Lemma}

\theoremstyle{definition}
\theoremstyle{remark}
\newcommand{\R}{\mathbb R}\newcommand{\Z}{\mathbb Z}
\newcommand{\T}{\mathbb T}\newcommand{\F}{\mathbb F}
\newcommand{\dist}{\operatorname{dist}}\newcommand{\covol}{\operatorname{covol}}
\newcommand{\adj}{\operatorname{adj}}
\newcommand{\wt}{\operatorname{wt}}\newcommand{\rank}{\operatorname{rank}}
\newcommand{\norm}[1]{\left\lVert #1\right\rVert}

\newcommand{\file}[1]{\nolinkurl{#1}}
\newcommand{\doi}[1]{\href{https://doi.org/#1}{\nolinkurl{#1}}}
\title{Euclidean View-Obstruction in Dimensions Eight to Ten\\
and the Asymptotic Deficit}
\author{}
\date{}
\makeatletter
\renewcommand{\maketitle}{%
  \begin{center}

    \vspace*{-1.0em}

    {\LARGE
      \@title
      \par
    }

    \vspace{3.15em}

    {\large
      Hu Tan\textsuperscript{1}
      \hspace{2.0em}
      Ying Zhang\textsuperscript{2}
      \par
    }

    \vspace{0.85em}

    \begin{minipage}[t]{0.49\textwidth}
      \centering
      \small
      \mbox{\textsuperscript{1}Academy of Mathematics and Systems Science}\\
      Chinese Academy of Sciences\\
      Beijing 100190, China\\[0.25em]
      \footnotesize
      \href{mailto:tanhu2020@amss.ac.cn}
      {\texttt{tanhu2020@amss.ac.cn}}
    \end{minipage}%
    \hfill
    \begin{minipage}[t]{0.49\textwidth}
      \centering
      \small
      \textsuperscript{2}School of Mathematical Sciences\\
      Soochow University\\
      Suzhou 215006, China\\[0.25em]
      \footnotesize
      \href{mailto:yzhang@suda.edu.cn}
      {\texttt{yzhang@suda.edu.cn}}
    \end{minipage}

    \vspace{0.9em}

  \end{center}

  \vspace{0.2em}
}
\makeatother
\begin{document}
\maketitle
\begin{abstract}
The Euclidean view-obstruction problem asks for the critical radius of
balls centered at the half-integer lattice that meet every positive ray
from the origin. Writing $\Delta_n$ for its square, we determine
$\Delta_8=8/17$, $\Delta_9=1087/1980$ and $\Delta_{10}=7/11$ by
computer-assisted proofs using Kahlon's seven-dimensional theorem.
The upper bounds combine full-support subspace averaging, projected
lattice coverings and exact finite certificates; dimension ten also
requires a finite search through short integer relations. Explicit
directions attain all three values, including a ten-dimensional
direction previously listed by Gritsenko, Skoruppa and Zagier.
In arbitrary dimension, we express their theta-block invariant in
geometric normalization and derive the known bounds
$c\log(n+1)\le n/12-\Delta_n\le C\log^3(n+1)$, with repeated coordinates
allowed throughout. Thus $\Delta_n\sim n/12$, while the deficit tends
to infinity more slowly than every positive power of $n$. The
supplementary archive contains the executable verifiers, exact
certificates and archived verification records.
\end{abstract}


\section{Introduction}
\subsection{The geometric question}
Consider a periodic array of identical opaque balls, with one ball centered
at every point of
\[
 c_n+\Z^n,\qquad c_n=(1/2,\ldots,1/2).
\]
An observer is at the origin. A viewing direction $u\in\R_{>0}^n$ is blocked
if the ray $\{tu:t>0\}$ meets at least one of the balls; tangency counts as
intersection. The problem asks for the threshold radius at which all such
directions are blocked. The half-lattice translation avoids placing an
obstacle at the observer and makes the coordinate symmetries explicit.
Only positive directions are considered, but repeated coordinates are
allowed throughout.

For a precise definition, let $r_n$ be the infimum of the radii that block
every positive ray, and put $\Delta_n=r_n^2$. If the reference ball has
diameter one and its dilation factor is denoted by $\nu(n)$, then
\begin{equation}\label{eq:normalization}
 r_n=\sqrt{\Delta_n},\qquad \nu(n)=2r_n,\qquad \nu(n)^2=4\Delta_n.
\end{equation}
This conversion is essential when comparing the older literature with
squared-distance computations.

The geometry has an equivalent periodic formulation. For
$v=(v_1,\ldots,v_n)\in\Z_{>0}^n$, define
\begin{equation}\label{eq:definitions}
 g(t)=\dist(t,1/2+\Z)^2=(\{t\}-1/2)^2,\quad
 F_v(t)=\sum_{i=1}^n g(v_it),\quad d(v)=\min_{t\in\T}F_v(t),
\end{equation}
where $\T=\R/\Z$ and $\{t\}$ denotes the fractional part of $t$. Then
\begin{equation}\label{eq:integer-sup}
 \Delta_n=\sup_{v\in\Z_{>0}^n}d(v).
\end{equation}
The reduction, including the distinction between an infimum and actual
intersection at the critical radius, is proved in Section~\ref{sec:geometry}.
Multiplication of $v$ by a positive integer does not change $d(v)$.
One may therefore restrict the supremum to primitive vectors, although
the finite censuses deliberately include all their multiples.

The integral of $g$ is $1/12$. Consequently $\Delta_n\le n/12$.
Determining the exact constant requires a direction that remains as far
as possible from the shifted lattice over an entire period, together with
an upper bound valid for infinitely many directions. The second requirement
is the central obstacle to a finite computation.

\subsection{Historical and analytic context}
The problem originates in Cusick's papers of 1972--1974
\cite{Cusick72,Cusick73,Cusick74}, which formulate view obstruction
in the geometry of numbers. The replacement of balls by axis-parallel cubes
connects the problem with simultaneous Diophantine approximation and the
lonely runner problem \cite{Wills,CusickPomerance}. For a cube of side
$\alpha$ centered at a half-integer point, intersection means that, at some
time $t$, all coordinates satisfy
\[
 \dist(tv_i,\Z)\ge(1-\alpha)/2.
\]
For a ball the condition instead bounds the sum of the squared deviations
from half-integers. This coupling of coordinates leads naturally to
projected lattices and covering radii.

For spheres, Dumir and Hans-Gill determined the three-dimensional
constant in 1986 \cite{DH86}. Chen studied dimensions three and four
in 1994 \cite{Chen3,Chen4}; Dumir, Hans-Gill and Wilker treated
dimension four in 1994 and dimension five in 1996
\cite{DHW94,DHW96}. The general theory of rational
subspaces, critical parameters and related isolation phenomena was developed
in \cite{DHW93,DHWII}; Chen also studied the connection with Kronecker's
theorem \cite{ChenK}. Kahlon's preprint \cite{Kahlon} proves the values in
dimensions six and seven using exact lattice searches and transference
bounds, including Banaszczyk's theorem \cite{Banaszczyk}. Our use of that
preprint is confined to its seven-dimensional upper bound and to the
context of the method. The computations in dimensions eight, nine and ten are separate from
those of Kahlon and are supplied with this paper.

A second connection is to theta blocks. Gritsenko, Skoruppa and Zagier
\cite[Sections~5--6]{GSZ} study the same minimum of a sum of periodic
quadratics, with a factor-of-two change of normalization. Their
Theorem~5.3 already establishes logarithmic and cubed-logarithmic bounds
for the corresponding deficit. We give a geometric formulation of their asymptotic theorem and a
detailed derivation that includes repeated coordinates. The analytic
inputs are Littlewood's inequality
\cite{MPS,Konyagin} and the integer-coefficient cosine-polynomial
construction of Belov and Konyagin \cite{BK}. The experimental table in
\cite[Table~2]{GSZ} also supplies
the lower bounds $37/80$ and $89/168$ in dimensions eight and nine.
The corresponding values considered here are $8/17$ and $1087/1980$.
The direction
$(1,1,2,3,4,5,6,7,8,9)$ and its ten-dimensional lower bound
$7/11$ already appear in \cite[Table~2]{GSZ}, after normalization.
For this direction, the further finite-dimensional task is to prove
the matching global upper bound.

\subsection{Main results}
The finite-dimensional contribution consists of the exact constants in
dimensions eight through ten, the finite reductions proving their upper
bounds, and explicit directions attaining them. The computations are
organized into three components: plane classification and a bounded
direction census in dimension eight; an enlarged catalogue, fourth-moment
constraints and a bounded census in dimension nine; and a finite
relation-lattice search in dimension ten. Appendix~\ref{app:index}
relates these proof components to the supplementary archive.

\begin{theorem}[Exact constants]\label{thm:exact}
The squared critical radii in dimensions eight, nine and ten are
\[
 \Delta_8=\frac8{17},\qquad
 \Delta_9=\frac{1087}{1980},\qquad \Delta_{10}=\frac7{11}.
\]
The respective primitive attaining directions can be chosen as
\[
 \begin{split}
 v^{(8)}&=(1,2,3,4,5,7,9,11),\\
 v^{(9)}&=(2,3,5,7,8,11,13,15,18),\\
 v^{(10)}&=(1,1,2,3,4,5,6,7,8,9).
 \end{split}
\]
For each of these three dimensions, balls of radius $\sqrt{\Delta_n}$
itself block every positive real ray.
\end{theorem}
The theorem supplies one attaining direction in each dimension.
The repeated entry in $v^{(10)}$ is intentional; equality directions
need not have distinct coordinates.
The dependency order of the upper bounds is
\[
 \Delta_7\ \longrightarrow\ \Delta_8\ \longrightarrow\
 \Delta_9\ \longrightarrow\ \Delta_{10}.
\]
The dimension-ten proof also uses the published KZ inequalities of
Pendavingh and van Zwam \cite{PVZ}, with their positive error terms retained.

\begin{theorem}[Known asymptotic bounds in geometric normalization]
\label{thm:asymptotic}
Set $E_n=n/12-\Delta_n$. There are absolute constants $c,C>0$ such that
for every $n\ge1$,
\begin{equation}\label{eq:asymptotic}
 c\log(n+1)\le E_n\le C\log^3(n+1).
\end{equation}
In particular,
\[
 \Delta_n\sim\frac n{12},\qquad E_n\longrightarrow\infty,
 \qquad E_n=o(n^\varepsilon)\quad(\varepsilon>0).
\]
\end{theorem}
All logarithms are natural. The bounds in \eqref{eq:asymptotic} are the
bounds of \cite[Theorem~5.3]{GSZ} after normalization. Section~\ref{sec:asymptotic}
gives a derivation from the two stated analytic inputs for arbitrary
positive integer multisets. Determining the precise order of $E_n$
remains open.

\begin{table}[ht]
\centering\small
\caption{Squared critical radii and one attaining direction. Rows $8$--$10$
give the exact constants established in this paper.}
\label{tab:values}
\vskip 6pt
\begin{tabular}{clll}
\toprule
$n$&$\Delta_n$&An attaining direction&Source\\\midrule
1&$0$&$(1)$&Immediate\\
2&$1/20$&$(1,2)$&\cite{Cusick74}\\
3&$3/28$&$(1,2,3)$&\cite{DH86}\\
4&$11/60$&$(1,1,2,3)$&\cite{Chen4,DHW94}\\
5&$41/168$&$(2,3,5,7,9)$&\cite{DHW96}\\
6&$107/312$&$(2,3,5,7,9,12)$&\cite{Kahlon}\\
7&$73/180$&$(1,3,4,7,8,10,11)$&\cite{Kahlon}\\
8&$8/17$&$(1,2,3,4,5,7,9,11)$&Section~\ref{sec:eight}\\
9&$1087/1980$&$(2,3,5,7,8,11,13,15,18)$&Section~\ref{sec:nine}\\
10&$7/11$&$(1,1,2,3,4,5,6,7,8,9)$&Section~\ref{sec:ten}\\
\bottomrule
\end{tabular}
\end{table}

\subsection{Proof strategy}
The proof has three stages. First, the closure of a real ray on the
torus reduces the geometric problem to rational subspaces. An averaging
lemma bounds a full-support subspace of dimension $d$ by $(n-d)/12$,
and a finite quadratic calculation evaluates any fixed integer direction.

Second, projected lattices transfer the remaining upper-bound problem
to Korkin--Zolotarev heights. Covering estimates and plane bounds constrain
their product. In dimensions eight and nine this produces a finite norm
bound for every possible counterexample. Dimension ten additionally
uses bounded short relations and a finite sequence of Gram-matrix
extensions. Exact finite certificates close these reductions; their verification
and the recorded coverage checks are described in
Section~\ref{sec:repro} and Appendix~\ref{app:index}.

Finally, in arbitrary dimension a Fourier expansion connects the
distance function to cosine sums. The McGehee--Pigno--Smith inequality
gives the logarithmic lower bound for the deficit, including repeated
speeds. A nonnegative integer-coefficient construction of
Belov--Konyagin gives the cubed-logarithmic upper bound. This argument
is independent of the finite certificates.

\section{Geometric reductions and exact distance evaluation}
\label{sec:geometry}
This section relates positive real rays to rational subspaces and gives
an exact finite calculation for a fixed integer direction. These two
reductions separate the global upper-bound problem from the verification
of explicit lower bounds.

\subsection{Rational subspaces and real rays}
For a rational subspace $W\subseteq\R^n$, write
\[
 \beta(W)=\dist(W,c_n+\Z^n)^2.
\]
Call $W$ full-support if it is contained in no coordinate hyperplane.
A subspace containing a positive vector has this property; the converse
need not hold. We use full-support statements because they remain invariant
under independent changes of coordinate signs.

\begin{proposition}\label{prop:rational}
The identity \eqref{eq:integer-sup} holds for the critical radius defined
as an infimum. For a positive real vector $u$, let $W$ be the smallest
rational linear subspace containing $u$. Then
\[
 \inf_{t>0}\dist(tu,c_n+\Z^n)^2=\beta(W).
\]
\end{proposition}
\begin{proof}
The continuous version of Kronecker's theorem identifies the closure of
the positive orbit $tu\bmod\Z^n$ with $(W+\Z^n)/\Z^n$; see the
view-obstruction formulation in \cite{DHW93,Kahlon}. Every tail of this
positive orbit has the same closure. The squared distance to $c_n+\Z^n$
is continuous on the torus, so its infimum on the orbit is its minimum
on that subtorus. This minimum equals $\beta(W)$.

Since $W$ is rational and $W\cap\R_{>0}^n$ is a nonempty relatively open
set, it contains a positive rational vector, hence a positive integer
vector $v$. Thus $\beta(W)\le d(v)$. Integer directions are themselves
allowed real directions, giving equality of the two suprema. A radius
strictly above their square root blocks every ray, and a radius strictly
below it misses some ray. This proves the assertion about the infimum.
\end{proof}
At the critical radius, an orbit with an unattained distance infimum
requires separate care. For dimensions eight through ten the strict
subspace margins established below settle this point.

\subsection{A full-support averaging lemma}
\begin{lemma}\label{lem:mean}
If $W\subseteq\R^n$ is a full-support rational subspace of dimension $d$,
then $\beta(W)\le(n-d)/12$.
\end{lemma}
\begin{proof}
Let $B$ be an integral basis matrix for $W\cap\Z^n$, with rows $b_i$,
and let $P=B(B^TB)^{-1}B^T$. For uniform $u\in\T^d$, set
$X_i=\{b_i u\}-1/2$. Each $b_i$ is nonzero, so $X_i$ has mean zero
and second moment $1/12$. If $b_i,b_j$ are linearly independent, their
joint character map is onto $\T^2$, and $\mathbb E X_iX_j=0$.
Otherwise write $b_i=pa$, $b_j=qa$ for a primitive integral row $a$
and nonzero integers $p,q$. The Fourier series of the periodic sawtooth,
or direct integration over a common period, gives
\[
 \mathbb E X_iX_j=\frac{\gcd(p,q)^2}{12pq},\qquad
 P_{ij}=pq\,a(B^TB)^{-1}a^T.
\]
Their product is nonnegative, including when $pq<0$. Hence
$\mathbb E\norm{PX}^2\ge\operatorname{tr}(P)/12=d/12$, and
$\mathbb E\norm{(I-P)X}^2\le(n-d)/12$.
Since $Bu\in W$ and $X=Bu-z-c_n$ for an integer vector $z$,
$(I-P)X$ is the negative of the orthogonal residual of a half-integer
center. Every such squared residual is at least $\beta(W)$.
\end{proof}

\subsection{Odd-center witnesses and complete periods}
For $v\in\Z_{>0}^n$, put $A=\norm v^2$. Every odd integer vector
$z\in(2\Z+1)^n$ gives the valid upper bound
\begin{equation}\label{eq:odd}
 d(v)\le \dist(z/2,\R v)^2
 =\frac{A\norm z^2-(v\cdot z)^2}{4A}.
\end{equation}
This inequality needs no assumption that $z/2$ is the nearest center
at the orthogonal projection time. To establish a lower bound, however,
one must minimize over all possible times, not merely exhibit a witness.

On $[0,1/2]$, form the ordered union of $0$, $1/2$ and all crossing
times $k/v_i$ in that interval. Inside each resulting cell $[a,b]$ the
nearest odd center is fixed, say $z$, and
\begin{equation}\label{eq:cell}
 F_v(t)=At^2-(v\cdot z)t+\norm z^2/4.
\end{equation}
Its minimum on the closed cell is obtained by clipping
$t_0=(v\cdot z)/(2A)$ to $[a,b]$. The minimum of these cell minima is
$d(v)$, since $F_v(1-t)=F_v(t)$. All quantities are rational.
This is also an independent way to audit a fast nearest-center sweep,
whose projection values are useful as upper witnesses.

For the three displayed directions, the numbers of half-period cells are
$16$, $32$ and $14$. Complete clipped calculations give the minima in
Theorem~\ref{thm:exact}. The squared norms are respectively $306$, $990$ and $286$. The last
direction reproduces the lower bound in \cite[Table~2]{GSZ}. Equality times in $[0,1/2]$ are
\begin{equation}\label{eq:times}
 \begin{array}{c|l}
 8&8/51,\ 25/51\\
 9&49/495,\ 94/495,\ 104/495,\ 149/495\\
 10&9/22.
 \end{array}
\end{equation}
The supplementary verifier regenerates all cells, odd centers, clipped
minimizers and exact values using rational arithmetic. These calculations
establish the three lower bounds independently of the upper-bound
censuses and relation searches.

\section{Lattice coverings and finite product bounds}
\label{sec:lattice}
We next bound the distance to a rational subspace through covering
radii of projected lattices. The resulting inequalities reduce a
global upper bound to a finite product calculation. Background on
reduced bases and reciprocal lattices can be found in \cite{Cassels,LLS}.
For a rational $d$-dimensional subspace $W$, let
$D=\covol(W\cap\Z^n)^2$. The lattice
$\Gamma=P_{W^\perp}\Z^n$ has rank $k=n-d$ and squared covolume $1/D$.
Indeed its dual is $W^\perp\cap\Z^n$, and complementary primitive
subgroups of the unimodular lattice $\Z^n$ have equal covolumes.
Choose a Korkin--Zolotarev (KZ) basis of $\Gamma$, and let
$h_1,\ldots,h_k$ be its squared Gram--Schmidt lengths. Then
\begin{equation}\label{eq:kz-basic}
 \prod_{i=1}^k h_i=1/D,\qquad h_{i+1}\ge\tfrac34h_i.
\end{equation}
The second assertion follows from the shortest-vector property in
each projection and a size-reduction coefficient of absolute value
at most $1/2$.

\begin{lemma}[Fiber covering]\label{lem:fiber}
Let $W_r$ be obtained by adjoining integer lifts of the first $r$
projected basis vectors to $W$. If $W$ is full-support, then
\begin{equation}\label{eq:fiber}
 \beta(W)\le\beta(W_r)+\rho_r^2
 \le\frac{k-r}{12}+\rho_r^2,
\end{equation}
where $\rho_r$ is the covering radius of the lattice generated by
those first $r$ projected vectors in their span.
\end{lemma}
\begin{proof}
Choose a half-integer center minimizing the normal distance to $W_r$.
Adding an integer combination of the chosen lifts leaves this normal
component unchanged and changes the component in $W_r\cap W^\perp$
by an element of the prefix lattice. Choose such a combination with
remaining component of norm at most $\rho_r$. The two errors are
orthogonal. The second inequality is Lemma~\ref{lem:mean}.
\end{proof}

\begin{lemma}[Adjacent two-dimensional blocks]\label{lem:blocks}
For adjacent KZ heights $h,k$, the corresponding projected rank-two
lattice has covering radius squared at most
\begin{equation}\label{eq:block}
 \frac{h+k}{4}-\frac5{12}(h-k)_+.
\end{equation}
For a matching $M$ of disjoint adjacent pairs in $\{1,\ldots,r\}$,
define
\begin{equation}\label{eq:U}
 U_{r,M}(h)=\frac14\sum_{i=1}^r h_i
 -\frac5{12}\sum_{i\in M}(h_i-h_{i+1}).
\end{equation}
Then $\rho_r^2\le U_{r,M}(h)$.
\end{lemma}
\begin{proof}
Write the two-dimensional basis as $(\sqrt h,0)$ and
$(\mu\sqrt h,\sqrt k)$, with $0\le\mu\le1/2$ and
$k+\mu^2h\ge h$. Triangulating a fundamental parallelogram gives
the squared circumradius bound
\[
 \rho^2\le\frac{(k+\mu^2h)(k+(1-\mu)^2h)}{4k}.
\]
This expression decreases with $\mu$ on $[0,1/2]$: after scaling
$h=1$, its derivative has the sign of
$(2\mu-1)(k-\mu(1-\mu))$, which is nonpositive because $k\ge3/4$.
If $k\ge h$, the value at $\mu=0$ suffices. If $k<h$, set
$z=\sqrt{1-k/h}\le1/2$ and use $\mu\ge z$ to obtain
$\rho^2\le h/[2(1+z)]$. The required comparison reduces to
$z(1-2z)(3+2z)\ge0$.
Round successive blocks from the end of a prefix, adding the
orthogonal squared errors. Finally replace $(h_i-h_{i+1})_+$ by
$h_i-h_{i+1}$, which only enlarges the upper bound.
\end{proof}

\subsection{Why a finite vertex certificate is sufficient}
The inequalities \eqref{eq:kz-basic} and \eqref{eq:fiber}--\eqref{eq:U}
produce rational polyhedra in the positive heights. The objective is
a lower bound on their product, hence an upper bound on $D$ or $A$.
For the unbounded polyhedra used here, the following detail is needed.
Every coordinate has a positive lower bound, and every recession
direction $d$ satisfies $d_1\ge0$, $d_{i+1}\ge3d_i/4$. Thus all its
coordinates are nonnegative. By the vertex--recession decomposition
of a pointed polyhedron \cite{Ziegler}, any feasible point is a convex
combination of vertices plus a nonnegative recession vector. Increasing
coordinates cannot decrease their product. Concavity of
$\sum_i\log h_i$ then shows that the product is bounded below by
the smallest vertex product. This argument does not presume compactness.

The vertex-traversal programs first relax each right-hand side by a
fixed, tiny positive rational number. The feasible region expands,
so a lower bound on the expanded region remains valid for the original
one. Exact inverse-matrix and slack checks establish simplicity of
each visited vertex. Releasing each active constraint identifies
every outgoing edge by rational ratio comparisons, or certifies it
as unbounded. Closure of the queue, together with connectivity of
the vertex graph, establishes exhaustiveness. The perturbation is
not trusted generically: nondegeneracy is checked at every vertex.
All products and comparisons use arbitrary-precision integers.

\section{Dimension eight}\label{sec:eight}
Put $T_8=8/17$. We first reduce a classification of full-support planes
to finite checks, then bound the norm of a putative counterexample.
The supplementary archive contains the exact certificates for the
plane classification, the norm bound and the final direction census.

\begin{proposition}[Eight-dimensional plane certificate]\label{prop:8plane}
Every full-support rational plane $W\subset\R^8$ satisfies
$\beta(W)\le43/100$, unless, up to signed coordinate permutations,
its saturated integer lattice has basis matrix
\begin{equation}\label{eq:8plane}
 B=\begin{pmatrix}
 0&3\\1&-2\\1&-1\\1&1\\1&2\\2&-1\\2&1\\3&0
 \end{pmatrix}.
\end{equation}
For this exceptional plane, $B^TB=21I_2$, $D=441$, and
$\beta(W)\le19/42$.
\end{proposition}
\begin{proof}[Proof and finite certificate]
Suppose $\beta(W)>43/100$. There are six projected KZ heights.
The ordinary nearest-plane bound in \eqref{eq:fiber}, with
$x_i=3h_i$, implies
\[
 \sum_{i=1}^r x_i>r-21/25,\qquad x_{i+1}\ge3x_i/4.
\]
Relax the strict inequalities. At a tight prefix $j$, subtraction
of the preceding prefix inequality gives $x_j\le1$, whereas the
next prefix gives $x_{j+1}\ge1$. Thus the ratio constraint immediately
following a tight prefix cannot also be tight. At a vertex, active
ratio edges consequently form geometric blocks, each terminated by
a tight prefix. The blocks are precisely the $32$ compositions of
$6$; the first block has sum equal to its length minus $21/25$, and
each subsequent block has sum equal to its length. The positive
coordinate bounds and the recession argument of
Section~\ref{sec:lattice} reduce the product minimum to these vertices.
Their exact products give $D<6843$, hence $D\le6842$.

Choose a Gauss-reduced basis $a,b$ of $W\cap\Z^8$, and write
\[
 s=\norm a^2,\quad t=a\cdot b,\quad u=\norm b^2,
 \quad |2t|\le s\le u,\quad D=su-t^2.
\]
Then $3s^2\le4D$, so $s\le95$. If $s=1$, $W$ contains a
coordinate axis, and deletion of that coordinate, followed by
\eqref{eq:integer-sup} in dimension seven, gives
$\beta(W)\le\Delta_7<43/100$.

The remaining finite classification is implemented by
\file{verify_planes.js}. It first covers all $174\,047\,114$
sorted positive eight-tuples of squared norm at most $4000$, retaining
only $195$ not certified below $43/100$. Let $H$ be their norm set.
No norm $s\le95$ is in $H$, so $a$ cannot have full support.
For every $k$ with $\norm{b+ka}^2\le4000$ whose norm is not in $H$,
$b+ka$ must have a zero coordinate. Distinct $k$ require distinct
nonzero coordinates of $a$; there are at most $\min(7,s)$ of them.
All relevant $k$ lie in $[-46,46]$. This leaves ten Gram matrices.

For each forced zero at $k$, a row of the basis matrix is
$(q,-kq)$ with $q\ne0$. Enumerating these rows and all completions
with the prescribed Gram matrix yields $63$ row multisets, of which
$45$ are saturated. Saturation is checked by the gcd of all
$2\times2$ minors. Forty-four have explicit odd-center witnesses
at most $43/100$; the remaining matrix is \eqref{eq:8plane}.
For it, the all-ones odd vector gives
$[8-(11^2+3^2)/21]/4=19/42$. These finite row and witness checks
are retained in \file{planes_certificate.json}.
\end{proof}

Suppose now that a primitive $v$ has $d(v)>T_8$, and put $A=\norm v^2$.
The first projected KZ vector lifts to a plane $W\supset\R v$.
Since this vector is primitive in the projected lattice, the lift
together with $v$ is a saturated basis of $W\cap\Z^8$, so $D=Ah_1$.
If $W$ is exceptional, the one-dimensional fiber bound gives
\[
 T_8<d(v)\le19/42+441/(4A),\qquad A<157437/26<6056.
\]
Otherwise, $h_1>4(T_8-43/100)=69/425$. For every prefix and adjacent
matching, \eqref{eq:fiber} implies
\begin{equation}\label{eq:8height}
 \sum_{i\le r}h_i-\frac53\sum_{i\in M}(h_i-h_{i+1})
 \ge\frac{17r-23}{51}.
\end{equation}
Together with the six adjacent KZ constraints and the $h_1$ bound,
these are $60$ inequalities. The exact outward-perturbed vertex
check visits $2933$ vertices, $19\,830$ directed bounded edges and
$701$ unbounded edges, proving $\prod h_i>1/13700$.
Thus every strict counterexample must have $A\le13699$.

The complete direction census at this bound contains
$18\,766\,118\,836$ sorted positive eight-tuples. It certifies
\begin{equation}\label{eq:8certificate}
 17\bigl[A\norm z^2-(v\cdot z)^2\bigr]\le32A
\end{equation}
for odd-center witnesses, either individually or for whole families
by the extension bounds in Appendix~\ref{app:tails}. The recorded
disjoint counts are
\[
\begin{array}{lr}
\text{Deletion using the seven-dimensional theorem}&507\,096\,272\\
\text{Whole tails from a seven-coordinate center}&18\,246\,848\,326\\
\text{Adaptive whole-tail bounds}&10\,229\,385\\
\text{Individual additional centers}&1\,944\,853.
\end{array}
\]
Their sum agrees with an independent generating-function count;
the $790\,944\,734$ extendable seven-coordinate prefixes also agree.
No counterexample remains. The exact period evaluation of $v^{(8)}$
supplies the matching lower bound, proving $\Delta_8=T_8$.

\section{Dimension nine}\label{sec:nine}
Put $T_9=1087/1980$. Here a larger low-norm catalogue and fourth-moment
constraints make the plane classification finite and manageable.
The supplementary archive records the catalogue, the plane
certificates and the complete partitioned census used below.

\subsection{High directions and exceptional planes}
The complete census of sorted positive nine-tuples with $A\le10000$
contains $29\,044\,309\,143$ tuples. Retaining every tuple not proved
at most $51/100$ gives $2368$ directions and $1365$ distinct squared
norms. Each retained distance is checked by the clipped rational
cell method. Their maximum is $T_9$; within this catalogue its
attainers are $v^{(9)},2v^{(9)},3v^{(9)}$.
This local observation is not a global uniqueness assertion.

\begin{proposition}[Nine-dimensional plane certificate]\label{prop:9plane}
A full-support rational plane in $\R^9$ either has
$\beta(W)\le51/100$ or admits, after signed coordinate permutations,
one of the seven bases in Table~\ref{tab:9planes}.
In particular, every such plane has $\beta(W)\le31/60$.
\end{proposition}
\begin{table}[ht]
\centering\small
\renewcommand{\arraystretch}{1.65}
\caption{Exceptional nine-dimensional plane bases, written as ordered
lists of rows. The last column is an upper bound, not an asserted exact offset.}
\label{tab:9planes}
\vskip 6pt
\begin{tabular}{cclc}
\toprule
$(s,t,u)$&$D$&Rows $(a_i,b_i)$&$\beta(W)\le$\\\midrule
$(10,-5,10)$&75&\shortstack[l]{$(0,1),(0,1),(1,-2),(1,-1),(1,-1),$\\$(1,0),(1,0),(1,1),(2,-1)$}&$31/60$\\[3pt]
$(10,5,10)$&75&\shortstack[l]{$(0,1),(0,1),(1,-1),(1,0),(1,0),$\\$(1,1),(1,1),(1,2),(2,1)$}&$31/60$\\[3pt]
$(10,0,45)$&450&\shortstack[l]{$(0,1),(0,4),(1,-3),(1,-2),(1,-1),$\\$(1,1),(1,2),(1,3),(2,0)$}&$23/45$\\[3pt]
$(16,-8,40)$&576&\shortstack[l]{$(0,1),(0,4),(1,-3),(1,-2),(1,1),$\\$(1,2),(2,-2),(2,-1),(2,0)$}&$37/72$\\[3pt]
$(16,8,40)$&576&\shortstack[l]{$(0,1),(0,4),(1,-2),(1,-1),(1,2),$\\$(1,3),(2,0),(2,1),(2,2)$}&$37/72$\\[3pt]
$(56,-28,56)$&2352&\shortstack[l]{$(0,3),(0,4),(1,-2),(1,1),(2,-1),$\\$(3,-3),(3,0),(4,-4),(4,0)$}&$43/84$\\[3pt]
$(56,28,56)$&2352&\shortstack[l]{$(0,3),(0,4),(1,-1),(1,2),(2,1),$\\$(3,0),(3,3),(4,0),(4,4)$}&$43/84$\\
\bottomrule
\end{tabular}
\end{table}
\begin{proof}[Proof and finite certificate]
If $\beta(W)>51/100$, the seven projected heights satisfy
\[
 12\sum_{i\le r}h_i+20\sum_{i\in M}(h_{i+1}-h_i)
 \ge4r-28+48(51/100).
\]
With the six KZ chain constraints this gives $59$ inequalities.
Exact traversal of $3181$ vertices, $21\,480$ directed bounded
edges and $787$ unbounded edges proves $D<18226$.
Thus a reduced Gram matrix has $D\le18225$ and $s\le155$.
The case $s=1$ is excluded using $\Delta_8<51/100$.

For the $1\,735\,940$ remaining Gram triples, set
$q(k)=sk^2+2tk+u$. If $q(k)\le10000$ but its norm is absent
from the high catalogue, $b+ka$ must have a zero coordinate.
The relevant range is $|k|\le71$. The maximum number of such
zeros is $\min(8,s)$ if $s$ is absent from the catalogue, and
$\min(9,s)$ otherwise. Moreover the residual Gram matrix after
subtracting the forced matrices
$\left(\begin{smallmatrix}1&-k\\-k&k^2\end{smallmatrix}\right)$
must be positive semidefinite. These tests leave $713$ triples.

Let $H_4(q)$ be the set of fourth moments of retained vectors
of norm $q$. The polynomial
$P(k)=\sum_i(b_i+ka_i)^4$ must have degree at most four, and
except at the allowed zero-coordinate positions must belong to
$H_4(q(k))$. If the zero-coordinate allowance is $c$ and $m$
positions have empty $H_4$, then any $c-m+5$ nonempty positions
contain five genuine full-support positions. The program checks
that enough positions exist, enumerates every such five-point
selection and value choice, and interpolates exactly. Integer
coefficient, moment and compatibility tests leave $282$ polynomials
on $261$ Gram triples.

Writing $P(k)=\sum_{j=0}^4p_jk^j$, complete enumeration of the
second and fourth moments of $(a_i)$ and $(|b_i|)$, followed by
all pairings and signs, enforces
\[
 4\sum_i a_i b_i^3=p_1,\quad6\sum_i a_i^2b_i^2=p_2,
 \quad4\sum_i a_i^3b_i=p_3.
\]
The Cauchy, convexity and rearrangement bounds used to shorten
this search are necessary conditions only. There are $449$
row matrices, $119$ saturated. For an odd vector $z$, put
$Q=\norm z^2$, $P=a\cdot z$, $R=b\cdot z$. Its squared residual is
\begin{equation}\label{eq:plane-witness}
 \frac{DQ-uP^2+2tPR-sR^2}{4D}.
\end{equation}
Saved witnesses bound $112$ matrices by $51/100$; the seven
remaining matrices have the bounds in Table~\ref{tab:9planes},
using the all-minus-ones vector. Saturation, Gram data and all
witnesses are rechecked with arbitrary-precision integers.
\end{proof}

\subsection{Exceptional directions, a norm bound, and the full census}
If $v$ is primitive in an exceptional plane with offset upper bound
$b_W$, the one-dimensional fiber inequality gives
\begin{equation}\label{eq:exception-fiber}
 d(v)\le b_W+\frac D{4A}.
\end{equation}
Thus $d(v)>T_9$ forces the bounds
\[
\begin{array}{c|cccc}
D&75&450&576&2352\\\hline
A<&37125/64&2970&570240/139&4074840/257\\
A\le&580&2969&4102&15855.
\end{array}
\]
Write $v=ma+nb$. Positive definiteness gives
$m^2\le Au/D$, $n^2\le As/D$, so a finite integer rectangle
suffices. Saturation makes $\gcd(m,n)=1$ equivalent to
primitivity of $v$. Discard zero coordinates and sort the absolute
values. The exact check covers $1186$ Gram-labelled cases modulo overall sign and $76\,136$
cells, all at most $T_9$; for $D=450$ it also includes the safe
extra boundary $A=2970$.

For a strict counterexample outside these planes, the first
projected KZ height satisfies
\[
 h_1>4(T_9-51/100)=386/2475.
\]
The eight-height product is $1/A$. The prefix constraints are
\begin{equation}\label{eq:9height}
 12\sum_{i\le r}h_i+20\sum_{i\in M}(h_{i+1}-h_i)
 \ge48T_9+4r-32.
\end{equation}
There are $87$ prefix matchings, seven chain constraints and the
additional first-height constraint: $95$ in total. The exact
vertex check has $17\,812$ vertices, $138\,862$ directed bounded
edges and $3634$ unbounded edges. It bounds the reciprocal product
strictly below $45960$, so $A\le45959$.

The full census at this bound partitions by the first three
coordinates into $79\,075$ disjoint jobs. It covers
$21\,484\,571\,760\,845$ sorted positive nine-tuples and finds
none above $T_9$. The disjoint coverage categories are
\[
\begin{array}{lr}
\text{Seven-coordinate whole intervals}&18\,999\,514\,214\,339\\
\text{Eight-coordinate intervals}&2\,017\,779\,983\,325\\
\text{Deletion using the eight-dimensional bound}&16\,453\,487\,495\\
\text{Uniform single tails}&450\,702\,840\,824\\
\text{Adaptive tails}&104\,693\,510\\
\text{Individual directions}&16\,541\,352.
\end{array}
\]
Each job satisfies the sum identity, and independent generating-function
counts agree job by job. The census uses the proved two-tail and convex
interval inequalities of Appendix~\ref{app:tails}; logical coverage
does not mean separate numerical minimization for every tuple.
Together with the exact lower bound from $v^{(9)}$, this proves
$\Delta_9=T_9$.

\section{Dimension ten}\label{sec:ten}
Put $T=7/11$. In this dimension, a bounded direction census is
combined with a finite search on short integer relations. The
upper bound uses the complete branch coverage specified below, with
exact certificates and independent set audits in the supplementary archive.

\subsection{The base census and the plane bound}
The complete census of sorted positive ten-tuples with $A\le20000$
contains $3\,519\,937\,718\,862$ tuples, split into $19\,048$
completed jobs, and certifies $d(v)\le T$ throughout. A separate
complete catalogue at $A\le15000$ and threshold $3/5$ contains
$895\,707\,071\,005$ tuples; exactly $111$ directions on $96$
norms exceed $3/5$. Every retained value has a rational clipped-cell
audit. The stronger three-tail bound in Appendix~\ref{app:three}
is part of the census proof.

\begin{proposition}\label{prop:10plane}
Every full-support rational plane $W\subset\R^{10}$ satisfies
$\beta(W)\le3/5$.
\end{proposition}
\begin{proof}[Proof and finite certificate]
Assuming $\beta(W)>3/5$, the eight projected KZ heights and
the block inequalities yield $D<32806$, by an exact traversal
of $19\,124$ vertices. Thus $D\le32805$ and a reduced Gram
matrix satisfies $s\le209$. If $s=1$, delete the coordinate
axis and use $\Delta_9<3/5$. If $s=2$, the short vector has
two entries $\pm1$. At a fixed value of the other parameter,
these two coordinates can be adjusted with remaining squared
error at most $1/8$. The other eight coordinates give
\[
 \beta(W)\le\Delta_8+1/8=81/136<3/5.
\]
For $s\ge3$, there is no high-catalogue norm at most $209$, so
$a$ must have a zero coordinate. Each $k$ for which
$sk^2+2tk+u\le15000$ is not a high norm forces a distinct
nonzero coordinate of $a$ to vanish in $b+ka$. There can be
at most $\min(s,9)$ such $k$. All $4\,136\,388$ admissible
Gram triples violate this necessary condition, by at least four.
A direct enumeration and an independent discriminant-based reverse
enumeration agree. This excludes every possible plane.
\end{proof}
For a primitive strict counterexample $v$, the projected lattice
$\Gamma=P_{v^\perp}\Z^{10}$ has nine heights with product $1/A$.
The plane bound gives
\begin{equation}\label{eq:10h1}
 h_1>4(T-3/5)=8/55.
\end{equation}

\subsection{Additional KZ inputs and short relations}
Besides the classical inequalities, the following hold for every
valid consecutive index range:
\begin{align}
 3h_i&\le4h_{i+1},&2h_i&\le3h_{i+2},\notag\\
 25h_i+36h_{i+1}&\le48h_{i+2}+(40+7/10^6)h_{i+3},\label{eq:pvz}\\
 5h_i&\le2h_{i+3}+(8+3/10^4)h_{i+4},\notag\\
 4h_i+3h_{i+2}&\le4h_{i+3}+(8+1/20000)h_{i+4}.\notag
\end{align}
The last three are precisely the proved versions of
\cite[Theorems~6.1--6.2]{PVZ}. Their error-free analogues are not
assumed. Nonnegative rational linear combinations of these
inequalities and \eqref{eq:U} give $\rho_r^2\le C_rh_r$, with
\begin{equation}\label{eq:Cr}
\begin{split}
(C_1,\ldots,C_9)=\biggl(&\frac14,\frac12,\frac56,\frac{29}{24},
\frac{505003}{300000},\frac{265003}{120000},\\
&\frac{26165423}{9000000},\frac{13280261}{3600000},
\frac{4020502145513}{855000000000}\biggr).
\end{split}
\end{equation}
The saved multipliers are checked by direct rational vector addition;
no numerical linear-programming output is accepted without this check.

Reverse the dual basis of $\Gamma$ to obtain a basis
$a_1,\ldots,a_9$ of $\Lambda=v^\perp\cap\Z^{10}$. Successively
reduce each new vector to a nearest representative modulo its
predecessors. This operation corresponds in the primal basis to
adding earlier vectors to later ones, preserving every projected
KZ lattice and height. The relation prefix
$M_k=\Z a_1+\cdots+\Z a_k$ is primitive in $\Z^{10}$.
Write $G_k$ for its Gram matrix and $D_k=\det G_k$; set $D_0=1$.
Its squared Gram--Schmidt lengths satisfy
\begin{equation}\label{eq:dual}
 g_j=1/h_{10-j},\qquad D_k=\prod_{j\le k}g_j,
 \qquad A=\frac{D_k}{h_1\cdots h_{9-k}}.
\end{equation}

Primitivity makes the integral pairing map onto $\Z^k$.
Since $a_i\cdot\mathbf1\equiv\norm{a_i}^2\pmod2$, the minimum
squared normal offset of the half-integer centers is
\begin{equation}\label{eq:parity}
 R_k=\frac14\min_{\substack{x\in\Z^k\\x_i\equiv(G_k)_{ii}\ (2)}}
 x^TG_k^{-1}x.
\end{equation}
The complementary subspace contains $v$ and is full-support, so
$R_k\le k/12$ for $k\le7$ and $R_8\le3/5$. More importantly,
the fiber bound is
\begin{equation}\label{eq:relationfiber}
 d(v)\le R_k+U_{9-k,M}(h).
\end{equation}
A particular parity-compatible $x$ gives an upper bound on $R_k$
and suffices in this inequality. It does not certify a lower bound
on $R_k$. The final proof has no exclusions based on a claim that
the characteristic minimum exceeds the subspace mean; all recorded
``mean-exclusion'' lists are empty.

Since $T<d(v)\le C_9h_9$, one has $\norm{a_1}^2<C_9/T<8$.
Norm one is impossible because $v$ has full support. If
$m=\norm{a_1}^2\in\{2,\ldots,7\}$, then
\[
 \norm{a_2}^2<\frac{C_8}{T-1/(4m)}+\frac m4<8.
\]
Thus the first two relations both have squared norm at most seven.
For a general prefix, putting $r=9-k$ and using a covering-radius
upper bound $\rho(M_k)^2$ yields
\begin{equation}\label{eq:nextnorm}
 \norm{a_{k+1}}^2<\frac{C_r}{T-R_k}+\rho(M_k)^2.
\end{equation}
If a rational height certificate also gives $h_r\ge L>0$, then
$\norm{a_{k+1}}^2\le1/L+\rho(M_k)^2$.
The strict bound uses $\lceil b\rceil-1$ as the integer upper
limit; the weak bound uses $\lfloor b\rfloor$.

\subsection{Complete Gram extension and exact exclusions}
Let $w_i=a_i\cdot a_{k+1}$ and $q=\norm{a_{k+1}}^2$.
Nearest reduction implies
\begin{equation}\label{eq:nearest}
 2|x\cdot w|\le x^TG_kx\quad(x\in\Z^k),\qquad
 |w_i|\le\lfloor(G_k)_{ii}/2\rfloor.
\end{equation}
Enumerating the whole latter box and every allowed $q$ is therefore
complete; tests with $x\in\{-1,0,1\}^k$ only prune by necessary
conditions. Negating the new vector permits the first nonzero
entry of $w$ to be positive. No additional symmetry deletion is
needed for completeness. A further necessary dual-KZ test is
\begin{equation}\label{eq:dualkz}
 x^T\adj(G_k)x\ge D_{k-1}\quad(0\ne x\in\Z^k).
\end{equation}
If a violating vector exists, it lies in the finite box
\[
 |x_i|\le\left\lfloor\sqrt{\frac{D_{k-1}(G_k)_{ii}-1}{D_k}}\right\rfloor.
\]
This follows from $x_i^2\le(G_k)_{ii}x^TG_k^{-1}x$ and integrality
of the adjugate quadratic form. The square-root endpoints are
checked by integer comparisons.

Each surviving Gram prefix fixes the last $k$ primal heights.
Equations \eqref{eq:pvz}, \eqref{eq:10h1}, the mean-fiber bounds
and \eqref{eq:relationfiber} define a rational polyhedron in the
remaining $r$ heights. The following exclusions are used.
\begin{enumerate}
\item A rational Farkas certificate proves infeasibility by a
nonnegative combination of the original inequalities. A numerical
solver may propose multipliers, but exact arithmetic accepts them.
\item A product certificate proves $A\le20000$, reducing the branch
to the completed base census. All free heights have positive lower
bounds and finite upper bounds from the fixed tail, so the product
minimum is attained at a vertex. For small $r$, every $r$-element
active set is tested with exact Cramer arithmetic, including
degenerate vertices. Earlier layers also use complete box trees
and lower chords for logarithms, with rational series remainder
bounds. Both children of every split must be covered.
\item A necessary integral embedding condition fails. A realizable
Gram matrix must equal $\sum_{j=1}^{10}c_jc_j^T$ for integral
columns $c_j$. Full support requires
$c_j^TG_k^{-1}c_j<1$ for every $j$, and primitivity requires the
gcd of all maximal minors to equal one. Complete finite searches
with positive-semidefinite residuals certify negative instances.
Timeouts remain unresolved, never excluded.
\end{enumerate}
The covering bounds for the relation lattice may use one-dimensional
rounding or exact adjacent Gauss-reduced blocks: a planar Gram
$\left(\begin{smallmatrix}p&t\\t&q\end{smallmatrix}\right)$ with
$p\le q$, $2|t|\le p$ has squared covering radius
$pq(p+q-2|t|)/(4(pq-t^2))$. Taking the minimum over all partitions
into singleton and adjacent pair blocks is valid even though the
dual basis need not itself be KZ-reduced.

Local embedding obstructions provide an additional inexpensive
filter. A primitive $k$-row integer matrix has row rank $k$ over
every $\F_p$, so
\begin{equation}\label{eq:localrank}
 \rank(G_k\bmod p)\ge2k-10.
\end{equation}
If equality holds at an odd prime, the determinant square class
of the nonsingular quotient must be $(-1)^{10-k}$.
Indeed the row space $U$ then contains $U^\perp$ as its radical
of dimension $10-k$. Splitting off that many hyperbolic planes
in the ambient determinant-one quadratic space leaves
$U/U^\perp$ with the stated determinant class. These are
necessary conditions, not sufficient embedding criteria.

\subsection{The special family and the final branch ledger}
The exact product computations through rank five give, for the
remaining branches,
\begin{equation}\label{eq:80216}
 A\le\left\lfloor4492125/56\right\rfloor=80216.
\end{equation}
One rank-six Gram matrix has exactly $24$ primitive full-support
embeddings in the finite classification. They are all of
$G_2\oplus B_2$ type. Signed coordinate permutations and Weyl
symmetries reduce their full-support integer directions to
\[
 (a,b,a+b,2a+b,3a+b,3a+2b,x,y,x+y,x-y),
 \quad a,b\ge1,\quad x>y\ge1,
\]
with norm
$A=8(3a^2+3ab+b^2)+3(x^2+y^2)$.
All $15\,093\,137$ parameter tuples with $A\le80216$ have an
odd-center witness satisfying
\[
 11\bigl[A\norm z^2-(v\cdot z)^2\bigr]\le28A.
\]
The family closes the whole Gram branch because all its embeddings
were classified, not just one representative.

The initial complete tables have $63$ rank-two Gram matrices and
$1943$ rank-three candidates. After exclusions, $290$ rank-three
parents generate $21\,839$ rank-four candidates. A safe superset
of $794$ rank-four parents generates $164\,085$ rank-five
candidates. Exact certificate and embedding checks leave $1717$
rank-five parents for the final stages.

\begin{table}[ht]
\centering\small
\caption{Complete high-rank relation ledger in dimension ten. Each row
is a disjoint classification of the candidate set at that rank.}
\label{tab:ledger}
\vskip 6pt
\begin{tabular}{lrrr}
\toprule
&Rank 6&Rank 7&Rank 8\\\midrule
Parents&$1717$&$2203$&$61$\\
Candidates passing dual KZ&$1\,126\,981$&$3\,970\,866$&$295\,391$\\
Local embedding exclusions&$2032$&$113\,648$&$37\,572$\\
Height certificates&$1\,088\,664$&$3\,811\,705$&$257\,819$\\
Exact product bounds&$34\,055$&$45\,415$&$0$\\
Complete embedding negatives&$26$&$37$&$0$\\
Complete special-family closure&$1$&$0$&$0$\\
Remaining for next rank&$2203$&$61$&$0$\\\bottomrule
\end{tabular}
\end{table}
The archived set audits reconcile the actual sets, in addition to
these cardinalities.
Rank seven uses parts $0,1,2,3,4,5,6$; rank eight uses parts
$5,40,41,42,43,44,45$. Their parent sets are disjoint and have
unions equal to the preceding full survivor selections. Each
excluded Gram is linked to its exact certificate, and every
survivor is passed to the next layer. The final gate and independent set audit both record zero unresolved
branches.
Any actual counterexample would produce a prefix in every one
of these complete tables, contradicting the empty rank-eight
remainder. The exact period calculation for $v^{(10)}$ supplies the matching lower
bound, proving $\Delta_{10}=7/11$.

\subsection{Intersection at the critical radius}
For an irrational positive direction, the rational orbit closure
from Proposition~\ref{prop:rational} has dimension at least two.
The established plane bounds and Lemma~\ref{lem:mean} give
\[
\begin{array}{c|cc|c}
n&\dim W=2&\dim W\ge3&\Delta_n\\\hline
8&19/42&5/12&8/17\\
9&31/60&1/2&1087/1980\\
10&3/5&7/12&7/11.
\end{array}
\]
Both subspace bounds are strictly smaller than $\Delta_n$.
Density of the positive orbit therefore gives an actual
intersection with the interior of a critical-radius ball.
For rational directions the minimum occurs in a period and can
be shifted to a positive time. Every orbit tail is dense in the
same subtorus, so intersections can also be taken arbitrarily far
along the ray. If only positive half-integer centers are retained
in the geometric formulation, this last observation ensures that
the same critical radii apply. This completes
Theorem~\ref{thm:exact} for real rays, including its boundary assertion.
\section{The deficit in general dimension}\label{sec:asymptotic}
For a positive integer multiset $v$, define
\[
 e(v)=\frac n{12}-d(v),\qquad E_n=\inf_{|v|=n}e(v).
\]
The averaging bound makes these quantities nonnegative. The results
in this section do not depend on any finite-dimensional census.

\subsection{Normalization of the theta-block invariant}
The function $B$ used in \cite{GSZ} is $g/2$. Accordingly its
extremal quantity $\max(n)$ is $\Delta_n/2$, and its definition
(5.1) gives the exact identity
\begin{equation}\label{eq:weight}
 \wt(n)=n/2-12\max(n)=6E_n.
\end{equation}
Here $\wt(n)$ is the extremal quantity as defined there, allowing
fractional powers of the eta function; no additional integral-weight
rounding is imposed. This is why their logarithmic bounds transfer
directly to the sphere problem. Their experimental theta-block
tables are useful constructions but are not substitutes for the
global upper bounds of Sections~\ref{sec:eight}--\ref{sec:ten}.

\subsection{Fourier expansion and a uniform comparison}
Put
\[
 \beta_2(t)=\{t\}^2-\{t\}+1/6=g(t)-1/12,
 \quad G_v(t)=\sum_i\beta_2(v_it),
\]
and define $C_v(t)=\sum_i\cos(2\pi v_it)$ and
$K(v)=-\min_{\T}C_v$. All integrals on $\T$ use probability Haar
measure. Twice integrating by parts gives
\begin{equation}\label{eq:fourier}
 \beta_2(t)=\frac1{\pi^2}\sum_{k\ge1}\frac{\cos(2\pi kt)}{k^2},
 \qquad G_v(t)=\frac1{\pi^2}\sum_{k\ge1}\frac{C_v(kt)}{k^2}.
\end{equation}
Both series converge absolutely and uniformly. In particular,
\begin{equation}\label{eq:transfer}
 e(v)\le K(v)/6.
\end{equation}
Writing $\mu$ for the M\"obius function, inversion also gives
\begin{equation}\label{eq:mobius}
 C_v(t)=\pi^2\sum_{k\ge1}\frac{\mu(k)}{k^2}G_v(kt).
\end{equation}
Substitution of \eqref{eq:fourier} proves this formula: the
coefficient of $C_v(rt)$ is
$r^{-2}\sum_{k\mid r}\mu(k)$, which vanishes except at $r=1$.
Absolute convergence justifies regrouping the double sum.

\pagebreak[3]
\begin{proposition}\label{prop:21}
For every multiset $v\in\Z_{>0}^n$,
$\norm{C_v}_{L^1(\T)}\le21e(v)$.
\end{proposition}
\begin{proof}
Write $e=e(v)$ and $H=G_v+e\ge0$, so $\int H=e$.
Separate the positive and negative M\"obius coefficients:
\[
 P(t)=\pi^2\sum_{\mu(k)=1}\frac{H(kt)}{k^2},\qquad
 N(t)=\pi^2\sum_{\mu(k)=-1}\frac{H(kt)}{k^2}.
\]
The Euler products give
\[
 \sum_{k\ge1}\frac{\mu(k)}{k^2}=\frac6{\pi^2},\qquad
 \sum_{k\ge1}\frac{|\mu(k)|}{k^2}
 =\frac{\zeta(2)}{\zeta(4)}=\frac{15}{\pi^2}.
\]
Thus $C_v=P-N-6e$, while $\int P=21e/2$ and
$\int N=9e/2$. Multiplication by an integer preserves Haar
measure, which justifies these integrals. Since $(C_v)_+\le P$
and $\int C_v=0$,
$\norm{C_v}_1=2\int(C_v)_+\le21e$.
\end{proof}
This convenient comparison supplies all that is needed below.
Its role is to transfer an $L^1$ lower bound to the geometric deficit.

\subsection{The lower bound, including repeated speeds}
We use the following form of the McGehee--Pigno--Smith inequality
\cite[Theorem~2]{MPS}: there is an absolute $b>0$ such that for
distinct integers $q_1<\cdots<q_N$ and complex coefficients $a_j$,
\begin{equation}\label{eq:MPS}
 \int_{\T}\left|\sum_{j=1}^N a_je^{2\pi iq_jt}\right|dt
 \ge b\sum_{j=1}^N\frac{|a_j|}{j}.
\end{equation}
Konyagin independently resolved Littlewood's problem
\cite{Konyagin}; the weighted formulation \eqref{eq:MPS} is
the input used in this proof.

Suppose $v$ has $s$ distinct positive speeds, of multiplicities
$m_1,\ldots,m_s$ summing to $n$. The exponential form of $C_v$
has $2s$ distinct frequencies, with coefficients $m_j/2$.
Consequently
\[
 \norm{C_v}_1\ge\frac b2\sum_{j=1}^{2s}\frac1j
 \gg\log(2s+1).
\]
Moreover every Fourier coefficient is bounded in absolute value
by the $L^1$ norm, giving
\[
 \norm{C_v}_1\ge\tfrac12\max_j m_j\ge n/(2s).
\]
If $s\ge\sqrt n$, the logarithmic bound suffices; otherwise
the multiplicity bound is larger than $\sqrt n/2$ and also
dominates a fixed positive multiple of $\log(n+1)$.
Thus uniformly over all multisets,
$\norm{C_v}_1\gg\log(n+1)$. Proposition~\ref{prop:21}, followed
by the infimum over $v$, gives $E_n\gg\log(n+1)$.

\subsection{The upper bound with prescribed total multiplicity}
The needed consequence of Belov--Konyagin
\cite[Theorem~0.5(2) and Corollary~5.4]{BK} is that,
for every integer $n\ge3$, there are nonnegative integers
$a_1\ge\cdots\ge a_n$, with $\sum_{j=1}^n a_j=n$, and a real
$a_0\ll\log^3n$ such that
\begin{equation}\label{eq:BK}
 a_0+\sum_{j=1}^n a_j\cos(jx)\ge0\qquad(x\in\R).
\end{equation}
This is the bound for $K_{\Z}^{\downarrow}(n)$, whose coefficient
sum is exactly $n$; a bound merely along a subsequence of lengths
would not suffice. Choose speed $j$ with multiplicity $a_j$,
omitting zero multiplicities. This gives a vector of exactly
$n$ coordinates and $C_v\ge-a_0$. Equation~\eqref{eq:fourier}
therefore gives $G_v\ge-a_0/6$, so
\[
 E_n\le e(v)\le a_0/6\ll\log^3n.
\]
Enlarging the absolute constants handles $n=1,2$ and proves
Theorem~\ref{thm:asymptotic}. The deep analytic inputs
\eqref{eq:MPS} and \eqref{eq:BK} are invoked as external theorems;
their original proofs are not reproduced here.

\subsection{Consequences and the limits of the bounds}
The exact identity
\begin{equation}\label{eq:radiusgap}
 \sqrt{n/12}-r_n=
 \frac{E_n}{\sqrt{n/12}+\sqrt{\Delta_n}}
\end{equation}
and $\Delta_n\sim n/12$ imply, for sufficiently large $n$,
\[
 c'\frac{\log n}{\sqrt n}\le\sqrt{n/12}-r_n
 \le C'\frac{\log^3n}{\sqrt n}\longrightarrow0,
 \qquad r_n/\sqrt n\longrightarrow1/\sqrt{12}.
\]
In particular, $E_n=o(\sqrt n)$, which excludes a positive
square-root term in the deficit.

The exact deficits from the three computations are
\[
 E_8=10/51,\qquad E_9=199/990,\qquad E_{10}=13/66.
\]
Their variation is compatible with the general asymptotic bounds and
does not by itself determine the growth rate of $E_n$.

There is also a structural limit to the $L^1$ argument. Put
$A(n)=\inf_{|v|=n}\norm{C_v}_1$. The previous proof gives
$A(n)\gg\log(n+1)$. Conversely, $v=(1,\ldots,n)$ gives
$C_v=(D_n-1)/2$, where
$D_n(t)=\sin((2n+1)\pi t)/\sin(\pi t)$ is the Dirichlet kernel.
For $0<|t|\le1/2$,
$|D_n(t)|\le\min(2n+1,1/(2|t|))$, whence
\[
 \norm{D_n}_1\le1+\log(2n+1),\qquad
 A(n)\le1+\tfrac12\log(2n+1).
\]
Thus $A(n)\asymp\log(n+1)$, and improving only the constant
in $E_n\ge cA(n)$ cannot yield a higher-order deficit lower bound.
Nor does a lower bound for a restricted monotone-coefficient
cosine problem transfer backwards through the one-sided inequality
\eqref{eq:transfer}. The precise order of $E_n$ remains beyond
the conclusions proved here.

\section{Verification and further questions}\label{sec:repro}
The computer-assisted proof combines mathematical reductions with exact
finite calculations. This section describes the role of the supplementary
archive (which is available at https://zenodo.org/records/24037416) 
and distinguishes its verification levels; 
the logical map is given in Appendix~\ref{app:index}. 

\subsection{Exact certificates and verification levels}
The archive contains executable sources, exact witnesses, rational
certificates, partition records and verification logs. Its main entry
point uses Python's standard library. The archived implementations also
use JavaScript, Node.js and checked C\# arithmetic. Python integers and
rational fractions, JavaScript arbitrary-precision integers and integer
cross-multiplication implement the decisive comparisons. In the JavaScript
censuses that use ordinary numbers, integer intermediate values are bounded
below $2^{53}$. Floating-point calculations may propose an endpoint or a
certificate, but acceptance uses an exact check.

The verification levels have distinct mathematical scopes. An integrity
check validates the file manifest, the source-archive hash and nested
archive checksums. The quick check also recomputes every cell for the
three lower-bound directions. The selected replay checks saved ledgers,
the nine-dimensional two-tail identities and plane witnesses, and the
ten-dimensional gate, short-relation bounds and three-tail identities.
The terminal replay additionally recomputes all rank-eight local and
height exclusions, using the archived candidate sets. Full
eight-dimensional verification regenerates its plane classification,
norm bound and direction census. Complete regeneration in dimensions
nine and ten additionally requires the large censuses and relation
generation procedures described in the archive's reproduction guide.

These levels should be interpreted through the corresponding proof
obligations. File identity establishes which evidence is being checked.
Population counts check the sizes of finite search domains. The distance
witnesses, exact exclusions and parent--child set reconciliations establish
the geometric bounds and the completeness of the finite search. Neither
hashes nor cardinalities alone provide those conclusions.

\subsection{Archived verification records}
The archive contains records of a complete eight-dimensional rerun,
selected nine-dimensional arithmetic and ledger audits, and
ten-dimensional independent candidate-set audits and terminal replays.
In particular, its terminal replay records all $37\,572$ local and
$257\,819$ height exclusions at rank eight. That replay verifies
the exclusion arithmetic for the archived candidate population; the
separate set-generation audits establish the population to which it
applies. The full nine-dimensional census and ten-dimensional base
censuses and relation generations have their own archived records.

The present revision aligns the mathematical exposition with those
sources and records. The large computations were not regenerated during
this revision. The analytic argument in Section~\ref{sec:asymptotic}
is independent of the finite censuses. Kahlon's seven-dimensional theorem
and the inequalities of Pendavingh and van Zwam remain external inputs.

\subsection{Further questions}
The finite method raises two further problems: classifying all attaining
directions and constructing smaller independent verifiers for the large
relation tables. In general dimension, the main analytic question is
to determine the order of $E_n$ between the logarithmic and
cubed-logarithmic bounds. A sharper estimate would need information
beyond improving only the constant in the $L^1$ comparison.

\section*{AI assistance declaration}
OpenAI GPT-6 Astra was used in the determination of the exact 
view-obstraction constants in dimensions eight to ten. 
OpenAI Codex assisted with the preparation of the manuscript, including
expository editing, literature checks, local exact-arithmetic checks
and the comparison of the manuscript with the supplementary archive.
The authors are responsible for the mathematical statements, proofs
and computational evidence.

\section*{Data and code availability}
The source code and exact verification data accompany this paper in a
separate archive (https://zenodo.org/records/24037416). 
It includes the verification entry point, reproduction
instructions, the original dimensional packages and the later audit
records. The mathematical paper is supplied separately as a standalone
LaTeX source; the code archive contains no LaTeX source files. Original
computational evidence is preserved byte for byte inside the archive.

\appendix
\section{Single-tail, two-tail and interval certificates}\label{app:tails}
\subsection{A single appended coordinate}
For a positive prefix $w$ and an odd center $z$, write
$S=\norm w^2$, $D=w\cdot z$, $Q=\norm z^2$, $N=SQ-D^2$.
Appending a coordinate $M$ and odd center coordinate $q$ gives
\[
 \dist((z,q)/2,\R(w,M))^2
 =\frac{N}{4S}+
 \frac{S}{4(S+M^2)}\left(q-\frac{MD}{S}\right)^2.
\]
The nearest odd integer to $MD/S$ is at distance at most one.
Hence
\begin{equation}\label{eq:onetail}
 d(w,M)\le\frac{N}{4S}+\frac{S}{4(S+M^2)}.
\end{equation}
Taking an optimal center for $w$ yields the deletion bound
$d(w,M)\le\Delta_{n-1}+S/[4(S+M^2)]$.
For $n=8$, $T=8/17$ and $\Delta_7=73/180$, a strict
counterexample requires $199M^2<566S$; retaining equality in a
census is a harmless enlargement of its candidate set.

\subsection{Two appended coordinates}
\begin{lemma}\label{lem:twotail}
Let $S=\norm w^2>0$, $L=\max_i w_i$, and $L\le x\le y$.
Then
\[
 d(w,x,y)\le d(w)+M(S/L^2),\qquad
 M(s)=\max\left(\frac14,\frac{(s+1)^2}{2(s+2)^2}\right).
\]
An arbitrary odd-center upper bound $N/(4S)$ may replace $d(w)$.
\end{lemma}
\begin{proof}
The normal prefix residual and the tail fiber are orthogonal.
The latter has Gram matrix
\[
 \frac1{S+x^2+y^2}
 \begin{pmatrix}S+y^2&-xy\\-xy&S+x^2\end{pmatrix}.
\]
Normalize $s=S/x^2$, $t=y/x\ge1$. If $t\ge s+1$, nearest-plane
rounding in a suitable order gives
\[
 \rho^2\le\frac14\left(\frac{s+1}{s+1+t^2}+\frac{s}{s+1}\right)
 \le\frac14\left(\frac1{s+2}+\frac{s}{s+1}\right)<\frac14.
\]
For $1\le t\le s+1$, the standard two-triangle tiling gives
\[
 \rho^2\le f(s,t)=
 \frac{(s+1)(s+t^2)(2s+(t-1)^2)}{4s(s+1+t^2)^2}.
\]
Put $u=t-1\in[0,s]$. If $s\ge\sqrt2$, the numerator of
$f(s,1)-f(s,t)$ after multiplication by the positive factor
$4s(s+1+t^2)^2(s+2)^2/(s+1)$ is
\[
 u\{4s^2(s+2)+(s^3+7s^2-4)u+(6s^2-8)u^2+(s^2-2s-4)u^3\}.
\]
If the final coefficient is negative, combine its term with the
linear term inside the braces, using $u^2\le s^2$. The result
is at least $u(s^4-s^3+3s^2-4)$, which is positive for
$s\ge\sqrt2$: its derivative is $s(4s^2-3s+6)>0$, and its
value at $\sqrt2$ is $6-2\sqrt2$. The other terms are nonnegative.

If $0<s\le\sqrt2$, the numerator of $1/4-f(s,t)$ after
multiplication by $4s(s+1+t^2)^2$ is
\[
 s(2-s^2)+4su+(-s^2+4s-1)u^2+(2s-2)u^3-u^4.
\]
For $s\le1$, the last three terms are at least $(2s-1)u^2$.
This is nonnegative when $s\ge1/2$; for smaller $s$ use
$4su-u^2\ge3su\ge0$. For $1\le s\le\sqrt2$, the cubic
term is nonnegative and the quadratic and quartic terms are
at least $(-2s^2+4s-1)u^2\ge0$. Therefore $f(s,t)\le M(s)$.
Finally, $M$ is nondecreasing and $S/x^2\le S/L^2$.
\end{proof}

\subsection{Exact whole-interval skipping}
For the nine-dimensional census, another valid two-tail bound is
the successive single-tail bound
\[
 \frac{N}{4S}+\frac{S}{4(S+L^2)}+
 \frac{S+L^2}{4(S+2L^2)}.
\]
With $T=1087/1980$ and $K=2^{20}$, define
\begin{align*}
 r_{\rm old}&=\lfloor4TK\rfloor-
 \left\lceil\frac{KS}{S+L^2}\right\rceil-
 \left\lceil\frac{K(S+L^2)}{S+2L^2}\right\rceil,\\
 r_{\rm new}&=\lfloor4TK\rfloor-
 \max\left(K,\left\lceil\frac{2K(S+L^2)^2}{(S+2L^2)^2}\right\rceil\right),
 \qquad r=\max(r_{\rm old},r_{\rm new}).
\end{align*}
Then $KN\le rS$ certifies every allowed two-tail extension.
When the last prefix coordinate increases, with preceding norm
$S_0$ fixed, $S/L^2=1+S_0/L^2$ decreases and both tail bounds
improve. Freezing the current $r$ is therefore conservative.

For a fixed odd center $(z_0,q)$ and variable last prefix
coordinate $y$, its projection numerator is
\[
 N(y)=N_0+S_0q^2-2D_0qy+Q_0y^2.
\]
The inequality $KN(y)\le r(S_0+y^2)$ has a convex quadratic
left-minus-right side, since $KQ_0-r>0$ in the implemented
seven-coordinate case ($Q_0\ge6$, $r<3K$). Exact verification
at two endpoints proves the entire interval. Eight-coordinate
intervals use the same reasoning with the single-tail bound.
Every proposed endpoint is checked by integer arithmetic; a
floating-point quadratic root is only a proposal. If no witness
succeeds, enumeration continues down to the complete period test.

\section{The three-tail covering certificate}\label{app:three}
The ten-dimensional census uses the following finite polynomial
certificate for all real tail parameters, not a numerical grid.

\begin{lemma}\label{lem:threetail}
Let $S>0$, $0<L\le x\le y\le z$ and $0<S/L^2\le7$.
For $a=(x,y,z)$ and $A=S+x^2+y^2+z^2$, the lattice with Gram
$G=I_3-aa^T/A$ has covering radius satisfying
\[
 \rho(G)^2\le H(S/L^2),\qquad
 H(s)=\max\left(\frac12,\frac{3s^2+9s+8}{4(s+3)^2}\right).
\]
\end{lemma}
\begin{proof}[Proof with explicit polynomial verification]
Normalize the parameters by setting $s=S/x^2$, $t=y/x$ and
$u=z/x$. Thus the normalized parameters are $s,1,t,u$.
Since $s\le S/L^2$ and $H$ is nondecreasing, it suffices to prove
the bound for $0<s\le7$, $1\le t\le u$. In the following
calculation write $z$ for the normalized variable $u$.

If $z(1+t)>s+1+t^2$, first quotient by the third basis vector
and then round in its direction. The squared covering radius is
at most
\[
 \rho_2^2(s;1,t)+\frac{s+1+t^2}{4(s+1+t^2+z^2)}.
\]
For $t\ge s+1$ this is less than $1/4+1/6<1/2$.
Otherwise use $f(s,t)$ from Lemma~\ref{lem:twotail} and bound
the second term by $(1+t)^2/[4(s+1+t^2+(1+t)^2)]$.
Let $C_0=s+1+t^2$, $R=(1+t)^2$. Positivity of the difference
from $1/2$ reduces to positivity of
\[
 P=2sC_0^2(C_0+R)
 -(s+1)(s+t^2)(2s+(t-1)^2)(C_0+R)-sC_0^2R.
\]
With $s=7u$, $t=1+7uv$, the polynomial $P/u$ has bidegree
$(6,6)$ on $[0,1]^2$. All $49$ tensor-Bernstein coefficients
are positive; their minimum is $140$. The exact verifier
constructs the polynomial, computes these coefficients and
re-expands them to recover the original identity.

In the complementary region $z(1+t)\le s+1+t^2$, the entire
domain has the parametrization
\[
 t=1+sv,\qquad z=t+\frac{s(1-v)w}{1+t},\qquad 0\le v,w\le1.
\]
The six coordinate-order path simplices triangulate a fundamental
cube. The maximum of their circumradii bounds the lattice
covering radius: for barycentric coordinates $\lambda_j$ in
a simplex of circumcenter $q$ and radius $R_\pi$,
\[
 \sum_j\lambda_j\norm{p-v_j}_G^2=R_\pi^2-\norm{p-q}_G^2\le R_\pi^2.
\]
Reversal of the path gives the same radius, so three orders
$\pi=(1,2,3),(1,3,2),(2,1,3)$ suffice. Put
\[
 H_0=1+t,\quad b=(H_0,tH_0,tH_0+s(1-v)w),
 \quad S_0=sH_0^2,\quad A_0=S_0+\sum b_i^2,
\]
and in order $\pi$ set
$E_i=A_0-b_{\pi_i}(2\sum_{j<i}b_{\pi_j}+b_{\pi_i})$.
Solving the circumcenter equations with
$G^{-1}=I+bb^T/S_0$ gives
\[
 R_\pi^2=\frac{N}{4D},\qquad
 N=S_0\sum_i E_i^2+\left(\sum_i b_{\pi_i}E_i\right)^2,
 \quad D=S_0A_0^2.
\]
For $0<s\le5$, substitute $s=5u$. The inequality
$R_\pi^2\le1/2$ is $2D-N\ge0$; after removing a factor $u$,
each of the three polynomials has tridegree $(10,10,4)$ and
all $605$ tensor-Bernstein coefficients are nonnegative.
For $5\le s\le7$, substitute $s=5+2u$. The inequality becomes
\[
 (3s^2+9s+8)D-(s+3)^2N\ge0.
\]
Each polynomial has tridegree $(13,10,4)$ and all $770$
tensor-Bernstein coefficients are nonnegative. In total the
obtuse-region check has $3(605+770)=4125$ coefficients, with
independent monomial re-expansion for every identity.

For clarity, these checks use the exact basis-conversion identity
\[
 x^k=\sum_{i=k}^d\frac{\binom{i}{k}}{\binom{d}{k}}
 \binom di x^i(1-x)^{d-i}
\]
in each variable. All arithmetic is rational, and every Bernstein
basis polynomial is nonnegative on the unit interval. Thus the
$4125+49=4174$ coefficient signs prove the continuum inequalities.
\end{proof}
For a sorted seven-coordinate prefix, $S/L^2\le7$ automatically.
The normal-prefix residual is orthogonal to the three-tail fiber,
so $d(w,x,y,z)\le N/(4S)+H(S/L^2)$. With $X=L^2$, $K=2^{20}$,
the integer certificate
\[
 KN\le S\left[\lfloor4TK\rfloor-
 \max\left(2K,\left\lceil
 \frac{K(3S^2+9SX+8X^2)}{(S+3X)^2}\right\rceil\right)\right]
\]
therefore proves an entire three-coordinate tail safe. This
is the envelope required by the dimension-ten census. The polynomial
identities and coefficient signs above establish its continuum
validity once checked in exact arithmetic.

\section{Verification map for the supplementary archive}\label{app:index}
The archive follows the same dependency order as the proof. Its
eight-, nine- and ten-dimensional packages contain the original sources
and exact data. A later ten-dimensional review package contains
independent candidate-set generation and exclusion checks. The following
map describes the evidence needed at each stage without tying the
mathematical argument to individual implementation file names.

\begin{longtable}{@{}>{\raggedright\arraybackslash}p{0.22\textwidth}>{\raggedright\arraybackslash}p{0.73\textwidth}@{}}
\toprule
Proof component&Verification evidence and acceptance condition\\\midrule\endhead
Candidate lower bounds&Complete rational half-period cells for the three
displayed directions, with clipped minimizers and the values
$8/17$, $1087/1980$ and $7/11$.\\[5pt]
Dimension eight&Exact plane witnesses, the height-product certificate
giving $A\le13699$, and a complete census of $18\,766\,118\,836$
sorted positive directions. The seven-prefix count is $790\,944\,734$;
no counterexample remains.\\[5pt]
Dimension nine&The high-direction catalogue, fourth-moment and plane
certificates, exceptional-plane checks, and the height-product bound
$A\le45959$. The final census has $79\,075$ partitions and covers
$21\,484\,571\,760\,845$ directions, with zero counterexamples.\\[5pt]
Dimension-ten base&The census at $A\le20000$ covers
$3\,519\,937\,718\,862$ directions in $19\,048$ partitions. The
high catalogue at $A\le15000$ and the reverse Gram enumeration support
the plane bound.\\[5pt]
Covering and norm bounds&Exact rational combinations of the KZ
inequalities, parity-compatible witnesses and the $4174$ Bernstein
coefficient checks in Appendix~\ref{app:three}.\\[5pt]
Relation generation&Complete bounded Gram extensions and parent--child
set audits. Ranks six, seven and eight contain respectively
$1\,126\,981$, $3\,970\,866$ and $295\,391$ candidate Grams after
the stated dual-KZ tests.\\[5pt]
Terminal exclusions&Every candidate is linked to an exact exclusion
certificate or to the next-rank survivor set. Rank seven uses parts
$0,1,2,3,4,5,6$; rank eight uses parts $5,40,41,42,43,44,45$.
The terminal split is $37\,572$ local and $257\,819$ height exclusions,
with zero unresolved branches.\\[5pt]
Special family&All $24$ primitive full-support embeddings are accounted
for, followed by the complete check of $15\,093\,137$ parameter tuples
at $A\le80216$.\\
\bottomrule
\end{longtable}

The executable guide specifies commands, parameters and expected outputs
for each component. The independent set audits compare the mathematical
sets passed between stages, while the saved logs document the original
completed computations. New runs use fresh working directories so that
the original evidence remains available for comparison.

\end{document}